\documentclass[english,russian,12pt,twoside]{article}
\usepackage{amssymb,amsmath}
\usepackage[russian]{babel}
\usepackage[cp866]{inputenc}
\begin{document}

\title{МЕТОД ОБОБЩЕННЫХ ФУНКЦИЙ В \\НЕСТАЦИОНАРНЫХ КРАЕВЫХ ЗАДАЧАХ \\ДЛЯ  ВОЛНОВОГО УРАВНЕНИЯ}
\author{\bf{Алексеева Людмила А.}}
\date{}
\maketitle \centerline{\textit {Институт математики МОН
РК,\,ул.Пушкина, 125,Алма-Ата,050100 Казахстан} }
\centerline{alexeeva@math.kz } \vspace{10mm}

\begin{abstract}
Рассматривается многомерный аналог уравнения Даламбера в
пространстве обобщенных функций. Излагается метод построения
условий на фронтах ударных волн. Разработан метод обобщенных
функций для построения решений начально-краевых задач для волновых
уравнений в пространствах разной размерности. Строятся динамические
аналоги формул Грина и Гаусса для решений волновых уравнений в
пространстве обобщенных функций. Построены их регулярные
интегральные представления и сингулярные граничные интегральные
уравнения для решения начально-краевых задач математической физики.
\end{abstract}

Решение многих задач акустики, гидромеханики, теории упругости и
других  разделов физики связано с решением краевых задач для
волнового уравнения - многомерного аналога уравнения Даламбера,
описывающего процессы распространения волн в однородных изотропных
средах, поэтому весьма  актуально   построение эффективных
способов их решения для областей с произвольной геометрией и
разнообразным видом граничных условий. Применение численных
методов конечных разностей, конечных элементов или других,
связанных с дискретизацией области решения,  приводит к
разрешающим системам высокого порядка. Определенные сложности
возникают при построении  разностных сеток и обеспечении точности
выполнения граничных условий и условий на фронтах ударных волн,
где  производные функций терпят разрыв.

Наиболее эффективным методом исследования таких задач является
метод граничных интегральных уравнений. Главное его достоинство -
снижение на порядок размерности решаемых уравнений, что весьма
существенно, особенно для областей, содержащих бесконечно
удаленные точки. Существующие методы и программы
сплайн-апроксимации произвольных контуров снимают проблему
ограничения формы рассматриваемых областей при использовании МГИУ.
В настоящее время  этот  метод широко используется для решения
стационарных  задач математической физики, что связано с успехами
в разработке теории ГИУ для эллиптических уравнений и систем.

Решение   нестационарных динамических задач на основе метода ГИУ
требует введения понятия обобщенного решения, что связано с
особенностью фундаментальных решений волновых уравнений, которые
принадлежат классу обобщенных функций. Кроме того классическое
понятие дифференцируемости решений для гиперболических уравнений
резко сужает класс задач, полезных для приложений. В частности,
типичные физические процессы, сопровождающиеся ударными волнами,
не описываются дифференцируемыми решениями гиперболических
уравнений.

Здесь  излагается метод ГИУ для построения решений
на\-чаль\-но-краевых задач для волновых уравнений в пространствах
разной размерности.  Строятся динамические аналоги формул Грина
для решений краевых задач для волнового уравнения в пространстве
обобщенных функций. Построены их регулярные интегральные
представления и сингулярные граничные интегральные уравнения для
решения соответствующих на\-чаль\-но- краевых задач.

\textbf{1. Обобщенные решения волнового уравнения,  ударные
волны.} Рассматривается многомерный  аналог уравнения Даламбера

\begin{equation}\label{(1.1)}
\Box_c u \equiv \Delta u - \frac{1}{c^ 2} \frac{{\partial ^2
u}}{{\partial t^2 }} = G(x,t),\quad x \in R^N ,\,\;t \in R^1 .
\end{equation}
Здесь $\Box_c $ - волновой оператор (даламбертиан), $\Delta $-
оператор Лапласа, $G$-- локально интегрируемая функция.

Уравнение (\ref{(1.1)}) строго гиперболическое, класс его решений
содержит разрывные по производным функции. Поверхности разрыва $F$
в $R^{N + 1}$--это характеристические поверхности уравнения
(\ref{(1.1)}), которые удовлетворяют характеристическому уравнению
в пространстве $R^{N + 1}  = \{(x,\tau  \equiv ct )\}$:
\begin{equation}\label{(1.2)}
\nu _\tau ^2  - \sum\limits_{j = 1}^N {\nu _{\rm j}^2 }  = 0,
\end{equation}
где $\nu (x,\tau ) = \left( {\nu _1 ,...,\nu _N, \nu _\tau  }
\right)$ - вектор нормали к $F$. Ему соответствует конус
характеристических нормалей - световой конус, для которого  $\nu
_\tau   = \nu _{N + 1}  < 0$ [1,2]. В $R^N $ такие поверхности
движутся  с единичной скоростью по $\tau$ :
\begin{equation}\label{(1.3)}
1 =  - \nu _\tau  /\left\| \nu  \right\|_N ,\quad \left\| \nu
\right\|_N  = \sqrt {\nu _j \nu _j }
\end{equation}
(по повторяющимся индексам $i,j$ в произведении здесь и далее
всюду проводится суммирование от 1 до $N$). В  $R^N $ им
соответствуют волновые фронты  $F_t$, движущиеся  со скоростью $c$
по времени  $t$ . На них выполняются условия непрерывности Адамара
:
\begin{equation}\label{(1.4)}
\left[ {u\left( {x,t} \right)} \right]_{F_t }  = 0,\quad \left[
{\dot u + cn_i \,u,_i } \right]_{F_t }  = 0,
\end{equation}
где через $\left[ {f\left( {x,t} \right)} \right]_{F_t } $
обозначен  скачок  $f$  на $F_t$ :
   \[
\left[ {f\left( {x,t} \right)} \right]_{F_t }  = f^ +  \left(
{x,t} \right) - f^ -  \left( {x,t} \right) = \mathop {\lim
}\limits_{\varepsilon  \to  + 0} \left( {f\left( {x + \varepsilon
n,t} \right) - f\left( {x - \varepsilon n,t} \right)}
\right),\quad x \in   F_ t,
\]
$ n(x,t)$ - единичный вектор нормали к $ F_ t $, направленный в
сторону распространения фронта волны:
\begin{equation}\label{(1.5)}
n_i  = \frac{{\nu _i }}{{\left\| \nu  \right\|_N }} =
\frac{{grad\,F_t }}{{\left\| {grad\,F_t } \right\|}},\;i =
1,...,N;\;
\end{equation}
Последнее равенство справедливо, если уравнение фронта волны можно
представить в виде $ F_t (x,t) = 0 $ при условии существования
$grad\,F_t $.

Класс подобных решений гиперболических уравнений называют
\textit{ударными волнами}, на их фронтах производные функций и
даже сами функции могут терпеть скачки.

 Из второго условия (\ref{(1.4)})  следует, на фронтах
\begin{equation}\label{(1.4*)}
\dot u^ -   + cn_i \,u,_i^ -   = \dot u^ +   + cn_i \,u,_i^ +
\end{equation}
Если перед фронтом волны  $u \equiv 0 $ (среда в покое), это
равенство дает полезное соотношение на фронте волны:   $
(grad\,u,n) =  - c^{ - 1} \dot u,\quad x \in F_t . $ Заметим, что
касательные производные к характеристической поверхности, в силу
непрерывности $u$,  также  непрерывны, т.е.
\begin{equation}\label{(1.6)}
\left[ {u,_\tau  \gamma _\tau   + u_j \gamma _j } \right]_F  =
0\quad    \textrm{ для }\,\,\forall \gamma \in R^{N+1}:\,\,(\nu
,\gamma ) = 0.
\end{equation}
В частности, если $\gamma  = \gamma ^j  = ( - \nu _j ,\,\,\nu
_\tau \delta _1^j ,\,\,\nu _\tau  \delta _2^j ,\,\,\nu _\tau
\delta _2^j )$, это приводит к условиям вида:
\begin{equation}\label{(1.6*)}
\left[ { - u,_\tau  \nu _j  + u_j \nu _\tau  } \right]_F  = 0
\Rightarrow \;\left[ {\dot un_j  + cu_j } \right]_{F_t }  = 0
\end{equation}

Далее рассмотрим функции $u(x,t)$, которые непрерывны вместе с
производными до второго порядка включительно почти всюду, за
исключением конечного или счетного числа поверхностей разрыва -
волновых фронтов, достаточно гладких почти всюду, на которых
выполнены условия на скачки (1.3). Назовем такие решения
\textit{классическими}. Покажем, что они являются обобщенными
решениями уравнений (\ref{(1.1)}).

Для этого  рассмотрим   (\ref{(1.1)})  на прост\-ранст\-ве
обоб\-щен\-ных функ\-ций  $ D'(R^{N + 1} ) = \\=\left\{ {\hat
f(x,\tau )} \right\}$, определенных на пространстве $D(R^{N + 1} )
= \left\{ {\varphi (x,\tau )} \right\}$ бесконечно
дифференцируемых финитных функций [2]. Значение $ \hat f$  на $
\varphi $, как принято ,  обозначаем $(\hat f,\varphi )$. Для
регулярных
 функций $\hat f$, соответствующих локально
интегрируемым  $f$,\\ $ (\hat f,\varphi ) = \int\limits_{R^{N + 1}
} {f(x,\tau )} \varphi (x,\tau )dx_1 ...dx_N d\tau. $ Далее
обозначим $ dV(x) = dx_1 ...dx_N $.

О п р е д е л е н и е.  Функция $ \hat f \in D'(R^N )$  называется
\textit{обобщенным решением} уравнения (\ref{(1.1)}) , если для
любого $ \varphi  \in D(R^{N + 1} )$ выполняется равенство: $
(\Box_c \hat f,\varphi ) \equiv (\hat f,\Box_c \varphi ) =
(G,\varphi ).$

Л е м м а 1.1. \textit{Если $u$-- классическое решение
(\ref{(1.1)}) , то $ \hat u$ является его обобщенным решением.}

Д о к а з а т е л ь с т в о. Если $u$ имеет конечный разрыв на
$F$, то в $ D'(R^{N + 1} )$, согласно правилам  определения
обобщенной производной [2] ,
\[
\hat u,_j  = u,_j  + {\rm [}u{\rm ]}_F \nu _j \delta _F {\rm
(}x,\tau {\rm )},
\]
где первое слагаемое справа - классическая производная по   $x_j
,\;j = 1,...,N + 1,\quad \left\| {\bf \nu } \right\| = 1,\;\delta
_F $  - сингулярная обобщенная функция- \textit{простой слой} на
$F$:
\[
\left( {{\rm [}u{\rm ]}_F \nu _j \delta _F {\rm (}x,\tau {\rm
)}{\rm ,}\varphi {\rm (x}{\rm ,}\tau {\rm )}} \right) =
\int\limits_F {{\rm [}u{\rm (}x{\rm ,}\tau {\rm )]}_F \nu _j {\rm
(x}{\rm ,}\tau {\rm )}} \varphi {\rm (x}{\rm ,}\tau {\rm
)}dF(x,t), \quad\forall \varphi  \in D(R^{N + 1} ).
\]
Здесь  интеграл по $ F$ поверхностный.   В силу непрерывности $u$
вне фронта волны, $ \left[ u \right]_F  = \mathop {\lim
}\limits_{\varepsilon  \to  + 0} \left( {u(x + \varepsilon \,n,t)
- u(x - \varepsilon \,n,t)} \right) = u^ +  (x,t) - u^ -  (x,t) =
\left[ u \right]_{F_t }. $ Поэтому, с учетом (\ref{(1.5)}),
получим: $ \hat u,_j  = u,_j  + {\rm [}u{\rm ]}_{F_t } \nu _j
\delta _F {\rm (}x, \tau {\rm ) = }u,_j  + \left\| \nu  \right\|_N
{\rm [}u{\rm ]} _{F_t } n_j \delta _F , \,\, \hat u,_{jj}  =
u,_{jj}  + \left[ {u,_j } \right]_{F_t } \left\| \nu \right\|_N
n_j \delta _F  + \partial _j \left\{ {\left\| \nu  \right\|_N {\rm
[}u{\rm ]}_{F_t } n_j \delta _F } \right\}. $ В силу
(\ref{(1.3)}), $ \hat u,_\tau   = u,_\tau   + {\rm [}u{\rm ]}_{F_t
} \nu _\tau \delta _F  = c^{{\rm  - 1}}
 u,_t  - \left\| \nu  \right\|_N {\rm [}u{\rm ]}_{F_t }
 \delta _F $, $
\hat u,_{tt}  = \\c^{{\rm  - 2}} u,_{tt}  - c^{{\rm  - 1}} \left[
{u,_t } \right]_{F_t } \left\| \nu  \right\|_N \delta _F  -
\partial _\tau  \left\{ {\left\| \nu \right\|_N {\rm [}u{\rm
]}_{F_t } \delta _F } \right\}{\rm ,}\quad j = 1,...,N. $  С
учетом этих равенств и условий Адамара (\ref{(1.1)}) , получим
\[
\displaylines{ \Box  _c \hat u =\Box _c u + \left\{ {c^{{\rm  -
1}} \left[ {u,_t } \right]_{F_t }  + \left[ {n_j u,_j }
\right]_{F_t } } \right\}\left\| \nu  \right\|_N \delta _F {\rm
(}x,\tau {\rm )} +  \cr  c^{ - 1} \partial _t \left\{ {\left\| \nu
\right\|_N {\rm [}u{\rm ]}_{F_t } \delta _F {\rm (}x,\tau {\rm )}}
\right\} +
\partial _j \left\{ {\left\| \nu  \right\|_N {\rm [}u{\rm ]}_{F_t
} n_j \delta _F {\rm (}x,\tau {\rm )}} \right\} = \hat G(x,t), }
\]
поскольку все плотности простых и двойных слоев на $F_t $ равны
нулю. Действительно, второе слагаемое равно нулю в силу второго
условия (\ref{(1.4)}) на фронтах. А два других в силу первого,
поскольку их действие на   определяется  так: $$ \displaylines{
\left( {c^{ - 1} \partial _t \left\{ {\left\| \nu \right\|_N {\rm
[}u{\rm ]}_{F_t } \delta _F {\rm (}x,\tau {\rm )}} \right\} +
\partial _j \left\{ {\left\| \nu  \right\|_N {\rm [}u{\rm ]}_{F_t
} n_j \delta _F {\rm (}x,\tau {\rm )}} \right\},\varphi (x,t)}
\right) =  \cr = \left( {\left\| \nu \right\|_N {\rm [}u{\rm
]}_{F_t } \delta _F {\rm (}x,\tau ,\frac{{\partial \varphi
}}{{\partial n}} - \varphi ,_\tau  } \right). }$$ Что и
требовалось доказать.

\textit{Замечание 1.}  Из этой леммы следует, что условия на
фронтах ударных волн легко получить, рассматривая классические
решения гиперболических уравнений как обобщенные. Достаточно
приравнять нулю плотности соответствующих  независимых сингулярных
обобщенных функций - аналогов простых, двойных и др. слоев,
возникающих при обобщенном дифференцировании решений.  Определение
таких условий на основе классических методов весьма трудоемкая
процедура.

\textit{Замечание 2. } Уравнение (\ref{(1.1)})   допускает
обобщенные решения со скачком производных и на подвижных
поверхностях $F(x,t) = 0$, скорость движения которых может
зависеть от точки фронта: $v(x,t) =  - F,_t /\left\| {grad\,F}
\right\|$, если на них выполняются условия Адамара (\ref{(1.4)}) .
Такие решения могут порождаться  правой частью уравнения, если
носитель $G(x,t)$ расширяется с течением времени в $R^N $.

\textbf{2. Постановка краевых задач и единственность решений}.
Пусть в области $x \in S^ -  $, ограниченной поверхностью Ляпунова
S ([2], c. 409), строится решение (\ref{(1.1)}) при $t \ge 0$.
Обозначим $n = \left( {n_1 ,...,n_N } \right)$ - единичный вектор
внешней нормали к  $S$, $D^ -   = S^ -   \times R^ +  ,\;S^ - \in
R^N ,\;R^ +   = [0,\infty )$.

\textit{Начальные условия}:  при $t = 0$
\begin{equation}\label{(2.1)}
u(x,0) = u_0 (x)\,\,  \textrm{для} \,\, x \in S^ -   + S,  \quad
u,_t (x,0) = \dot u_0 (x) \textrm{ для} \,\, x \in S^- .
\end{equation}
Здесь рассмотрим  две краевые задачи, соответствующие условиям
Дирихле и Неймана.

\textit{Граничные условия}: \\(первая КЗ)
\begin{equation}\label{(2.2)}
 u(x,t) = u_S (x,t)\,\,  \textrm{для} \,\,x \in S,\;t \ge 0;
\end{equation}
(вторая КЗ)
\begin{equation}\label{(2.3)}
n_i (x)\partial _i u(x,t) = \frac{{\partial u}}{{\partial n}} =
p(x,t)\,\, \textrm{ для }\,\, x \in S,\;t \ge 0.
\end{equation}

Для первой краевой задачи  выполнены условия согласования
граничных и начальных данных:
\begin{equation}\label{(2.4)}
u_0 (x) = u_S (x,0)\,\, \textrm{для} \,\, x \in S.
\end{equation}

На волновых фронтах, если они возникают, выполняются условия
Адамара (\ref{(1.4)}) . Заметим, что ударные волны всегда
возникают, если не выполнено условие согласования начальных и
граничных данных по скоростям:
\begin{equation}\label{(2.1*)}
\dot u_0 (x) = \dot u_S (x,0)\,\,  \textrm{для }\,\, x \in S,
\end{equation}
что типично для физических задач (здесь и далее $\partial _t u =
\dot u$ ). В этом случае в начальный момент времени на границе $S$
формируется фронт ударной волны, который распространяется со
скоростью $c$ в $R^N$. Для построения непрерывно дифференцируемых
решений это условие является необходимым. Здесь мы его вводить не
будем. Предполагается, что начальные условия заданы и известно
одно из граничных условий соответственно рассматриваемой краевой
задаче.

Введем функции $\;E = 0,5\left( {u,_\tau ^2  + \sum\limits_{j =
1}^N {u,_j^2 } } \right)\;,\;L = 0,5\left( {u,_\tau ^2  -
\sum\limits_{j = 1}^N {u,_j^2 } } \right)\;$.

Л е м м а  2.1. \textit{Если $u$-- классическое решение
(\ref{(1.1)}) , то на фронтах}
\begin{equation}\label{(2.5)}
\left[ E \right]_{F_t }  =  - c^{ - 1} \left[ {\dot
u\frac{{\partial u}}{{\partial n}}} \right]_{F_t }
\end{equation}
\begin{equation}\label{(2.6)}
\left[ {L(x,t)} \right]_{F_t }  =c ^{ - 2} \left( {\dot u^ -   +
\frac{{\partial u^ -  }}{{\partial n}}} \right)\left[ {\;\dot u}
\right].
\end{equation}

Д о к а з а т е л ь с т в о. В силу равенства:  $[ab] = a^ +  [b]
+ b^ -  [a]$ , с учетом (\ref{(1.1)})  и (\ref{(1.6*)}), получим :
\[
 \left[ {cE + \dot u\frac{{\partial u}}{{\partial
n}}} \right] = \;\left[ {0,5\left( {c^{ - 1} \dot u^2  + cu,_j
u,_j } \right)\; + \dot u\frac{{\partial u}}{{\partial n}}}
\right] =\]
\[
\begin{gathered}
= 0,5\left[ {\dot u\left( {c^{ - 1} \dot u\; + \frac{{\partial
u}}{{\partial n}}} \right)} \right] + 0,5\left[ cu,_j u,_j \right]
+ 0,5[\dot u\frac{\partial u}{\partial n} ] = \\ =  0,5c^{ - 1}
\left( {\dot u^ +  \left[ {\dot u + c\frac{{\partial u}}{{\partial
n}}} \right] + \left( {\dot u^ - + cu,_j^ -  n_j } \right)\left[
{\dot u} \right]} \right) + 0,5u,_j^ +  \left[ {\;cu,_j  + \dot
un_j } \right]+\cr  + 0,5\left( {cu,_j^ - + \dot u^ -  n_j }
\right)\left[ u,_j^?  \right] =  0,5c^{ - 1} \left[ \dot u
\right]\left( {\dot u^ {-}   + cu,_j^ -  n_j } \right) + \\ +
0,5\left[ {u,_j } \right]\left( {cu,_j^ -   + \dot u^ - n_j }
\right) =   0,5cu,_j^ -  [ u,_j  + c^{ - 1} n_j \dot {u} ] +
0,5c^{ - 1} \dot u^ -  [{cn_j u,_j  + \dot u} ] = 0
\end{gathered}
\]
(здесь $n$ - нормаль к фронту в $R^N $). Отсюда следует формула
(\ref{(2.5)})  леммы.

Далее, поскольку $[a^2 ] = \left( {a^ +   + a^ -  } \right)[a]$ ,
и в силу (\ref{(1.1)})  и (\ref{(1.6)}), получим (\ref{(2.6)}):
\[
\displaylines{ \;\left[ L \right] = 0,5\left[ {u,_\tau ^2  -
\sum\limits_{j = 1}^N {u,_j^2 } } \right] = \;0,5\left( {u,_\tau ^
+   + u,_\tau ^ -  } \right)\left[ {u,_\tau ^{} } \right] -
0,5\left( {u,_j^ +   + u,_j^ -  } \right)\left[ {\;u,_j } \right]
=  \cr = 0,5\left( {u,_\tau ^ +   + u,_\tau ^ -  } \right)\left[
{u,_\tau ^{} } \right] + 0,5\left( {u,_j^ +   + u,_j^ -  }
\right)n_j \left[ {\;u,_\tau  } \right] =  \cr = 0,5\left\{
{\left( {u,_\tau ^ +   + n_j u,_j^ +  } \right) + \left( {u,_\tau
^ -   + n_j u,_j^ -  } \right)} \right\}\left[ {\;u,_\tau  }
\right] = c^{ - 2} \left( {\dot u^ -   + \frac{{\partial u^ -
}}{{\partial n}}} \right)\left[ {\;\dot u} \right] .}
\]

\textit{Замечание} 1. Условие (\ref{(2.5)}) легко можно получить,
рассматривая соответствующее уравнение в  $D'(R^{N + 1} )$: $ \hat
E,_\tau   - (u,_\tau  u,_j ),_j  =  - u,_\tau G + \left\{ {\left[
E \right]\nu _\tau   - \left[ {u,_\tau  u,_j } \right]\nu _j }
\right\}\delta _F (x,\tau ) =$  \\ $=  - u,_\tau G - \left\| \nu
\right\|_N \left\{ {[ E ] + [ {u,_\tau  u,_j } ]n_j }
\right\}\delta_{F_{t}}(x,t). $ Следовательно, чтобы (\ref{(2.5)})
выполнялось в  $D'(R^{N + 1} )$, необходимо, чтобы $ \left[ E
\right]_F  + \left[ {u,_\tau u,_j } \right]_F n_j  = \left[ {E +
c^{ - 1} \dot u\frac{{\partial u}}{{\partial n}}} \right]_F  = 0$,
т.е. выполнялась формула (\ref{(2.5)}).

\textit{Замечание} 2. Eсли перед фронтом волны $u \equiv 0$, то, с
уче\-том (\ref{(1.4*)}), имеем $\left[ {L(x,t)} \right]_{F_t } =
0$, т.е. в этом случае функция $L$   непрерывна.

Т е о р е м а 2.1. \textit{Если $u(x,t)$ -- классическое решение
краевой задачи, то}
\[
\int\limits_{S^ -  } {(E(x,t)}  - E(x,0))dV(x) =  -
\int\limits_0^t {dt} \int\limits_{D^ -  } {G(x,t)u,_t dV(x)}  +
\int\limits_0^t {\int\limits_S {\left( {\dot u_S (x,t)p(x,t)}
\right)dS(x)dt} }
\]

Д о к а з а т е л ь с т в о.  Умножая (\ref{(1.1)})  на $u,_\tau $
в области дифференцируемости, после простых преобразований
получим:
\begin{equation}\label{(2.7)}
E,_\tau   - (u,_\tau  u,_j ),_j  =  - u,_\tau  G.
\end{equation}

А теперь проинтегрируем (\ref{(2.7)}) по $D^ -  $,  с учетом
разбиения области интегрирования волновыми фронтами $F_k $.
Заметим, что первые два слагаемые можно рассматривать как
дивергенцию соответствующего вектора в  пространстве $R^{N + 1} $,
которая в областях между фронтами непрерывна. Поэтому, используя
теорему  Остроградского-Гаусса в $R^{N + 1} $, получим
\[
\displaylines{ \int\limits_{D^ -  } {E,_\tau  } dV(x)d\tau  -
\int\limits_{D^ -  } {(u,_\tau  u,_j ),_j dV(x)d\tau }  +
\int\limits_{D^ -  } {u,_\tau  G(x,\tau )dV(x)d\tau }  =  \cr =
\int\limits_{D^ -  } {u,_\tau  G(x,\tau )dV(x)d\tau }  +
\int\limits_{S^ -  } {(E(x,\tau )}  - E(x,0))dV(x) -
\int\limits_0^\tau  {\int\limits_S {\left( {u,_\tau  u,_j n_j }
\right)dS(x)d\tau } }  +  \cr + \sum\limits_{F_k }
{\int\limits_{F_k } {\left[ {E\nu _\tau   - u,_\tau  u,_j \nu _j }
\right]_{F_k } dF_k (x,\tau } } ) = 0\quad  \cr}
\]
($dF_k (x,\tau ))$ - дифференциал  площади поверхности в
соответствующей точке волнового фронта).   В силу (\ref{(1.3)})  и
(\ref{(2.5)}),  $\left[ {E\nu _\tau   - u,_\tau  u,_j \nu _j }
\right]_{F_k }  =  - c^{-1}{{\left\| \nu  \right\|_N }}\left[ {cE
+ \dot u\frac{{\partial u}}{{\partial n}}} \right] = 0$. Поэтому
последний интеграл равен нулю. С учетом обозначений для граничных
функций, отсюда получаем формулу теоремы. Ее следствием является
следующая теорема.

Т е о р е м а 2.2. \textit{Если классическое решение первой
(второй) краевой задачи существует, то оно единственно. }

Д о к а з а т е л ь с т в о. В силу линейности задачи, достаточно
доказать единственность нулевого решения. Для него $G=0$,
начальные и соответствующие граничные условия нулевые. Тогда, как
легко видеть,  из теоремы 2.1 следует: $\int\limits_{S^ - }
{E(x,t)} dV(x) = 0$. Поскольку $E$ - неотрицательна, следовательно
$E \equiv 0 \Rightarrow u = const.$ Из начальных условий следует
$u \equiv 0.$

Т е о р е м а 2.3. \textit{Если  $u(x,t)$ -- классическое решение
краевой задачи, то}
\[
\displaylines{ \int\limits_{D^ -  } {L(x,t)} dV(x)dt =
\int\limits_{D^ -  } {uG(x,t)dV(x)dt}  - \int\limits_0^t
{\int\limits_S {\left( {u_S (x,t)p(x,t)} \right)dS(x)dt} }  +  \cr
+ c^{ - 2} \int\limits_{S^ -  } {\left( {u\dot u(x,t) - u_0 \dot
u_0 (x)} \right)dV(x)}  \cr}
\]

Д о к а з а т е л ь с т в о.  Умножая (\ref{(1.1)})  на $u$, после
простых преобразований получим:
\begin{equation}\label{(2.8)}
L(x,\tau ) + (uu,_j ),_j  - (uu,_\tau  ),_\tau   = G.
\end{equation}
Проинтегрируем (\ref{(2.8)}) по области $D^ -  $ с учетом
разбиения ее волновыми фронтами $F_k $. Аналогично, как в теореме
2.1, используя теорему  Остроградского-Гаусса, имеем
\[
\displaylines{ \int\limits_{D^ -  } {L\,} dV(x)d\tau  =
\int\limits_{D^ -  } {uGdV(x)d\tau }  + \int\limits_{D^ -  }
{\left( {(uu,_\tau  ),_\tau   - (uu,_j ),_j } \right)dV(x)d\tau }
=  \cr = \int\limits_{D^ -  } {uG\,dV(x)d\tau }  + c^{ - 1}
\int\limits_{S^ -  } {\left( {uu,_t  - u_0 \dot u_0 }
\right)dV(x)}  - \int\limits_0^\tau  {d\tau \int\limits_S {uu,_j
n_j dS(x)} }  +  \cr + \sum\limits_{F_k } {\int\limits_{F_k }
{\left[ {uu,_\tau  \nu _\tau   - uu,_j \nu _j } \right]_{F_k }
dF_k (x,\tau } } ).}
\]
В силу условий Адамара (\ref{(1.4)}),последнее слагаемое равно
нулю. Поэтому, с учетом условий на границе (\ref{(2.2)}),
(\ref{(2.3)}), получим формулу теоремы.

\textbf{3.  Динамический аналог формулы Грина в $D'(R^{N + 1} ).$
} Для построения решения КЗ перейдем в пространство $D'(R^{N + 1}
)$. Для этого введем характеристическую функцию области
определения решения  $H_D^ -  (x,t) \equiv H_S^ -  \left( x
\right)H(t)$, где $H_S^ -  \left( x \right)$-- характеристическая
функция множества $S^ -  $ равная 0,5 на его границе $S$, $H(t)$
-- функция Хевисайда, равная 0,5 при $t=0$. $H_D^ -  $ -
характеристическая функция пространственно-временного цилиндра $D^
-   = \left\{ {S^ -   \times R_{^ +  } } \right\}$.  Легко
показать, что
\begin{equation}\label{(3.1)}
\frac{{\partial H_D^ -  }}{{\partial x_j }} =  - n_j \delta _S
(x)H(t),\quad \frac{{\partial H_D^ -  }}{{\partial t}} =  - n_j
H_S^ -  (x)\delta (t).
\end{equation}

Введем обобщенные функции $\hat u(x,t) = u(x,t)H_D^ -  (x,t),\hat
G(x,t) = G(x,t)H_S^ -  \left( x \right)H(t)$,  где
$u(x,t)$--классическое решение КЗ,   и рассмотрим действие
волнового оператора на эту функцию. Поскольку $[u]_S  =  - u,$
выполняя обобщенное дифференцирование, подобно тому, как в п.1, и
используя (\ref{(1.1)})  в области дифференцируемости, получим
\begin{equation}\label{(3.2)}
\begin{gathered}
\Box_c \hat u\left( {x,t} \right) =  - \frac{{\partial
u}}{{\partial n}}\delta _S \left( x \right)H\left( t \right) -
H\left( t \right)\left( {un_j \delta _S \left( x \right)}
\right),_j -\\ - c^{ - 2} H_S^ -  \left( x \right)u_0 \left( x
\right)\dot \delta (t) - c^{ - 2} H_S^ -  \left( x \right)\dot u_0
\left( x \right)\delta (t) + \hat G,
\end {gathered}
\end{equation}
где $\beta \left( {x,t} \right)\delta _S (x)H(t)$ -- простой слой
на боковой поверхности цилиндра $\left\{ {S \times R^ +  }
\right\}$,  $\delta (t)$-- функция Дирака, $\frac{{\partial
u}}{{\partial n}} = u,_{i} n_i $ -- производная по нормали $n$ к
$S$. Заметим, что плотности простых и двойных слоев здесь
определяются граничными условиями, часть из которых,  в
зависимости от решаемой КЗ, известна, и заданными начальными
условиями.

Как известно из теории обобщенных функций [2], решением уравнения
(\ref{(3.2)}) является свертка фундаментального решения уравнения,
с правой частью уравнения (\ref{(3.1)}). В качестве такого возьмем
фундаментальное решение $\hat U\left( {x,t} \right)$,
удовлетворяющее условиям:
\begin{equation}\label{(3.3)}
\Box_c U = \delta (x)\delta (t),
\end{equation}
(\textit{условия излучения})
\begin{equation}\label{(3.4)}
 U = 0\quad \textrm{при} \,\,t < 0 \,\,\textrm{и}\,\, \left\| x \right\| > ct  .
\end{equation}
Назовем его \textit{ функцией Грина} уравнения (\ref{(1.1)}) .
Решение (\ref{(3.1)}) представимо в виде следующей свертки:
\begin{equation}\label{(3.5)}
\begin{gathered}
\hat{u}\left( {x,t} \right) = u\left( {x,t} \right)H_S^ - \left( x
\right)H\left( t \right) =  - \hat U\left( {x,t} \right)
* \frac{{\partial u}}{{\partial n}}\delta _S \left( x
\right)H\left( t \right) - \left( {\hat U * un_j \delta _S \left(
x \right)H\left( t \right)} \right),_j -\cr - c^{ - 2} \left(
{\hat U\mathop  * \limits_x H_S^ -  \left( x \right)u_0 \left( x
\right)} \right),_t  - c^{ - 2} \hat U\mathop  * \limits_x H_S^ -
\left( x \right)\dot u_0 \left( x \right) + \hat G * \hat U,\\
\end{gathered}
\end{equation}
где cимвол "$\mathop  * \limits_x $" означает, что свертка берется
только по $x$, поскольку при взятии свертки можно воспользоваться
свойством $\delta$-функции. Причем решение (\ref{(3.5)}))
единственно в классе функций, допускающих свертку с $U$.

Запишем (\ref{(3.5)})при нулевых начальных данных и $G = 0$,
используя свойство дифференцирования сверток:
\begin{equation}\label{(3.5*)}
\hat u\left( {x,t} \right) =  - \hat U\left( {x,t} \right) *
\frac{{\partial u}}{{\partial n}}\delta _S \left( x \right)H\left(
t \right) - \hat U,_j  * un_j \delta _S \left( x \right)H\left( t
\right).
\end{equation}

Формула выражает решение КЗ через граничные значения искомой
функции и ее нормальной производной и аналогична формуле Грина для
решений уравнения Лапласа [2]. Однако, в силу особенностей
фундаментальных решений гиперболических уравнений на фронте волны,
вид которых зависит от размерности пространства,  ее интегральное
представление дает расходящиеся интегралы (во втором слагаемом).
Для построения ее регулярного интегрального представления введем
функции -- первообразные по  $t$ :
\begin{equation}\label{(3.6)}
\hat W\left( {x,t} \right) = \hat U\left( {x,t} \right) * \delta
\left( x \right)H\left( t \right) = \hat U\left( {x,t}
\right)\mathop  * \limits_t H\left( t \right)\,\, \Rightarrow \,\,
\partial _t \hat W\left( {x,t} \right) = \hat U\left( {x,t} \right)
;
\end{equation}
\[
\hat H\left( {x,m,t} \right) = \frac{{\partial \hat W\left( {x,t}
\right)}}{{\partial x_i }}m_i  = \frac{{\partial \hat W\left(
{x,t} \right)}}{{\partial m}}.
\]
Заметим, что в силу ограниченности носителя по $t$ слева, обе
свертки существуют.  Легко видеть, что они также являются
решениями (\ref{(1.1)})  при $G = H(t)\delta (x)$  и $G =
H(t)\frac{{\partial \delta \left( {x,t} \right)}}{{\partial m}}$
соответственно. Справедлива следующая  теорема.

Т е о р е м а 3.1. \textit{В $D'(R_{N + 1} )$ решение КЗ
удовлетворяет уравнению}:
\begin{equation}\label{(3.7)}
\hat u\left( {x,t} \right) =  - \hat U\left( {x,t} \right) *
\frac{{\partial u}}{{\partial n}}\delta _S \left( x \right)H\left(
t \right) - \left( {\hat W,_j  * \dot un_j \delta _S \left( x
\right)H\left( t \right)} \right) -
\end{equation}
\[
- \left( {\hat W,_j \mathop  * \limits_x u_0 \left( x \right)n_j
\left( x \right)\delta _S \left( x \right) - c^{ - 2} \hat
U\mathop  * \limits_x H_S^ -  \left( x \right)\dot u_0 \left( x
\right) - c^{ - 2} \left( {\hat U\mathop  * \limits_x H_S^ -
\left( x \right)u_0 \left( x \right)} \right),_t  + \hat G * \hat
U} \right)
\]
При нулевых начальных данных и $G = 0$
\begin{equation}\label{(3.8)}
\hat u\left( {x,t} \right) =  - \hat U\left( {x,t} \right) *
\frac{{\partial u}}{{\partial n}}\delta _S \left( x \right)H\left(
t \right) - \hat W,_j  * \dot un_j \delta _S \left( x
\right)H\left( t \right)
\end{equation}

Д о к а з а т е л ь с т в о.  Легко показать, пользуясь
определением производной обобщенной функции и непрерывностью $u$,
что
\[
\left( {un_j \delta _S \left( x \right)H\left( t \right)}
\right),_t  = \dot u(x,t)n_j (x)\delta _S \left( x \right)H\left(
t \right) + u(x,0)n_j (x)\delta _S \left( x \right)\delta \left( t
\right).
\]
Используя это равенство и свойство дифференцирования свертки,
имеем
\[
\begin{array}{l}
\left( {\hat U * un_j \delta _S \left( x \right)H\left( t \right)}
 \right),_j  = \left( {\hat W,_t  * un_j \delta _S \left( x \right)H\left( t \right)}
  \right),_j  = \hat W,_j  * \left( {un_j \delta _S \left( x \right)H\left( t \right)}
   \right),_t  =  \\
= \hat W,_j  * \left( {un_j \delta _S \left( x \right)H\left( t
\right)}
 \right),_t  = \hat W,_j  * \dot u(x,t)n_j (x)\delta _S \left( x \right)H\left( t \right)
  + \hat W,_j  * u(x,0)n_j (x)\delta _S \left( x \right)\delta \left( t
  \right).
\end{array}
\]
Поскольку $ \hat W,_j  * u(x,0)n_j (x)\delta _S \left( x
\right)\delta \left( t \right) = \hat W,_j \mathop  * \limits_x
u_0 (x)n_j (x)\delta _S \left( x \right) $, подставляя эти
соотношения в (\ref{(3.5)}), получим формулу теоремы.

Из теоремы следует, что решение  задачи полностью определяется
граничными значениями нормальной производной функции $u\left(
{x,t} \right)$ и ее скорости $\dot u  $. По аналогии с
представлением решений уравнения Лапласа, эти формулы можно
назвать динамическим аналогом формулы Грина.

Формула (\ref{(3.7)}) обладает преимуществом в сравнении с
(\ref{(3.5)}), так как позволяет сразу перейти к ее интегральной
записи без регуляризации подынтегральных функций на фронтах, как
ранее было предложено в [3]. Для $x \in S$ формула (\ref{(3.7)})
дает, как покажем далее, граничное интегральное уравнение для
решения начально-краевой задачи, если для $x \in S$ известны
\begin{equation}\label{(3.9)}
\dot u\left( {x,t} \right) = v\left( {x,t}
\right)\quad\quad\textrm{(первая КЗ)};
\end{equation}
либо
\begin{equation}\label{(3.10)}
\frac{{\partial u}}{{\partial n}} = p\left( {x,t}
\right)\quad\quad\,\,\textrm{(вторая КЗ)}.
\end{equation}
Если известно одно из этих значений, то решая ГИУ на границе,
находим второе, после чего по формуле (\ref{(3.7)}) можно найти
$u(x,t)$.

\textbf{4. Динамический аналог формулы Гаусса.} Введем функции
\[
U\left( {x,y,t} \right) = \hat U\left( {x - y,t} \right) , W\left(
{x,y,t} \right) = \hat W\left( {x - y,t} \right) , H\left(
{x,y,m,t} \right) = \hat H\left( {x - y,m,t} \right),
\]
которые, в силу свойств симметрии волнового оператора и
$\delta$-функции, удовлетворяют следующим соотношениям симметрии
\begin{equation}\label{(4.1)}
\begin{gathered}
U\left( {x,y,t} \right) = U\left( {y,x,t} \right) ,\,\, W\left(
{x,y,t} \right) = W\left( {y,x,t} \right) ,\quad \frac{{\partial
W}}{{\partial x_j }} =  - \frac{{\partial W}}{{\partial y_j }},
\cr H\left( {x,y,m,t} \right) =  - H\left( {y,x,m,t} \right) = -
H\left( {x,y, - m,t} \right).
\end{gathered}
\end{equation}

Л е м м а  4.1. \textit{В $D'(R_{N + 1} )$ динамический аналог
формулы Гаусса имеет вид}:
\begin{equation}\label{(4.2)}
- \hat W,_i \mathop  * \limits_x n_i \left( x \right)\delta \left(
x \right) - c^{ - 2} \left( {\hat U\mathop  * \limits_x H_S^ -
\left( x \right)} \right),_t  = H_D^ -  \left( {x,t} \right)
\end{equation}

Д о к а з а т е л ь с т в о. Свернем обе части уравнения
(\ref{(3.3)}) для $U$ с $H_D^-(x,t)$  , используя свойства
дифференцирования сверток и  (\ref{(3.1)}):
\[
- U,_j *n_j (x)\delta _S \left( x \right)H(t) - c^{ - 2} \left(
{U,_t *H_S^ -  \left( x \right)\delta (t)} \right) = H_D^ -  (x,t)
\]
Выполняя свертку по $t$, с учетом (\ref{(3.6)}), получим формулу
леммы.

Формула (\ref{(4.2)}) является аналогом известной формулы Гаусса
для потенциала двойного слоя (с.406, [2]), которая дает
интегральное представление характеристической функции множества с
использованием фундаментального решения уравнения Лапласа. Формулу
Гаусса часто используют для вывода ГИУ краевых задач для
эллиптических уравнений и систем. Аналогично можно использовать
динамический аналог формулы Гаусса для вывода ГИУ гиперболических
уравнений. Здесь для построения ГИУ воспользуемся динамическим
аналогом формулы Грина.

Интегральная запись этих соотношений зависит от размерности
задачи. Далее дадим интегральное  представление формул теоремы
3.1. и леммы 4.1 для пространств размерности N=1,2,3, наиболее
характерных для физических задач.

\textbf{5. ГИУ для плоских краевых задач}.    В случае $N=2$ имеем
плоскую задачу. Рассмотрим вначале задачу с нулевыми начальными
условиями.
    Обозначим $dS(y)$ - дифференциал длины дуги кривой $S$ в
точке $y$, $S_t (x) = \left\{ {y \in S,\quad r = \left\| {y - x}
\right\| < ct} \right\}, $ $ \,S_t^ -  (x) = \left\{ y \in S^ - ,
r < ct \right\},\, r = \left\| {x - y} \right\|,\,dV(y) = dy_1
dy_2  $.

Т е о р е м а  5.1. \textit{При $N=2$ функция $\hat u\left( {x,t}
\right)$, удовлетворяющая нулевым начальным условиям $\left( {u_0
= 0,\;\dot u_0  = 0} \right)$, имеет следующее интегральное
представление}: (первая форма)
\begin{equation}\label{(5.1)}
\hat u = \frac{c}{{2\pi }}\int\limits_0^t {d\tau }
\int\limits_{S_\tau  ( x )} \left( {\frac{\partial u(y,t -
\tau)}{\partial n( y )} + \frac{1}{r}\frac{\partial r}{\partial n
(y)}\tau \dot u(y,t - \tau ) }\right )\frac{dS(y)}{\sqrt{ (c^2
\tau ^2 - r^2)}} ,
\end{equation}
\textit{которое, при перемене порядка интегрирования, имеет
вид}:\\
(вторая форма)
\[
\begin{gathered}
\hat u\left( {x,t} \right) = \frac{c}{{2\pi }}\int\limits_{S_t
\left( x \right)} {dS\left( y \right)} \int\limits_{{r
\mathord{\left/ {\vphantom {r c}} \right.
\kern-\nulldelimiterspace} c}}^t {\frac{\partial u\left( {y,t -
\tau } \right)}\partial n\left( y \right)\frac{d\tau
}{\sqrt{\left( {c^2 \tau ^2  - r^2 } \right)}}}+ \\ +
\frac{c}{{2\pi }}\int\limits_{S_t \left( x \right)}
{\frac{1}{r}\frac{{\partial r}}{{\partial n\left( y
\right)}}ds\left( y \right)} \int\limits_{{r \mathord{\left/
{\vphantom {r c}} \right. \kern-\nulldelimiterspace} c}}^t
{\frac{{\tau \dot u\left( {y,t - \tau } \right)d\tau
}}{{\sqrt{\left( {c^2 \tau ^2  - r^2 } \right)}} }}.
 \end{gathered}
\]
\textit{Для  $x\in S$ второй интеграл справа сингуляр\-ный,
берется в смыс\-ле глав\-но\-го зна\-че\-ния.}

Д о к а з а т е л ь с т в о. Для $N=2$ функция Грина имеет
следующий вид ([2], c. 206)
\begin{equation}\label{(5.2)}
\hat U\left( {x,t} \right) =  - \frac{{cH\left( {ct - R}
\right)}}{{2\pi {\sqrt{\left( {c^2 \tau ^2  - r^2 } \right)}} }}
,\quad R = \sqrt {x_1^2  + x_2^2 }.
\end{equation}
Вычисляя по формулам (\ref{(3.6)}), найдем $\hat W$ и $\hat H$
\begin{equation}\label{(5.3)}
\hat W\left( {x,t} \right) =  - \frac{{H\left( {ct - R}
\right)}}{{2\pi }}\ln \left( {\frac{{ct - \sqrt {c^2 t^2  - R^2 }
}}{R}} \right) ,\quad \hat H\left( {x,m,t} \right) =  -
\frac{{\,\,ct\,H\left( {ct - R} \right)}}{{2\pi \,\sqrt {c^2 t^2 -
R^2 } }}\frac{{x_j m_j }}{{R^2 }} .
\end{equation}
Если записать свертки в (\ref{(3.8)}) в интегральном виде с учетом
этих представлений, то получим соотношения  теоремы. Заметим, что
для $ x \notin S$ интегралы справа являются регулярными и потому
для таких $x$ справа и слева стоят регулярные функции. Докажем
справедливость (\ref{(5.1)}) и для $x \in S$ с учетом определения
$H_S^ - (x)$.

Обозначим $\epsilon$-окрестность точки $x$  ($\varepsilon  \ll ct,
t>0)$  через $ \textrm{Ш}_\varepsilon  (x) = \left\{ {y:\;r <
\varepsilon } \right\},$ $ S_\varepsilon ^ -  \left( x \right) =
S^ - - \textrm{Ш}_\varepsilon \left( x \right), S_\varepsilon ^ +
\left( x \right) = S^ +   - \textrm{Ш}_\varepsilon  \left( x
\right),\,O _\varepsilon \left( x \right) = \left\{ {y \in S:\;r
\le \varepsilon } \right\}, \,\,S_\varepsilon   = S -
O_\varepsilon,$ $\textrm{Ш}_\varepsilon ^ -   = S^ - \cap
\textrm{Ш}_\varepsilon,\,\,\textrm{Ш} _\varepsilon ^ +   = S^ +
\cap \textrm{Ш}_\varepsilon,\,\,\textrm{Г} _\varepsilon ^ \pm
\left( x \right) = \left\{ {y \in S^ \pm  :\;r = \varepsilon }
\right\}, r = \left\| {y - x} \right\|,\quad r,_j =
\frac{{\partial r}}{{\partial y_j }} = \frac{{y_j  - x_j }}{r},
\frac{{\partial r}}{{\partial n(y)}} = r,_j n_j (y).$

Пусть $x^*\in S$. Трансформируем контур $S$  в окрестности точки
$x^*$, обходя ее по $\varepsilon $-полуокружности в $S^{-}\,
(\varepsilon <<ct, t>0)$. Запишем динамический аналог формулы
Грина для контура $S_\varepsilon   +\textrm{Г} _\varepsilon ^ - $
в точке $x=x^*$ :
\[
0 = \frac{c}{\pi }\int\limits_{S_\varepsilon  \left( x^* \right) +
\textrm{Г}_\varepsilon ^ -  \left( x^* \right)} {H\left( {ct - r}
\right)dS\left( y \right)\int\limits_{{r \mathord{\left/
{\vphantom {r c}} \right. \kern-\nulldelimiterspace} c}}^t
{\frac{{\partial u\left( {y,t - \tau } \right)}}{{\partial n\left(
y \right)}}\frac{{d\tau }}{{\sqrt{c^2 \tau ^2  - r^2}} }}} +
\]
\[
+ \frac{c}{{2\pi }}\int\limits_{S_\varepsilon  \left( x^* \right)}
{H\left( {ct - r} \right)\frac{1}{r}\frac{{\partial r}}{{\partial
n\left( y \right)}}dS\left( y \right)\int\limits_{{r
\mathord{\left/ {\vphantom {r c}} \right.
\kern-\nulldelimiterspace} c}}^t {\frac{{\tau \dot u\left( {y,t -
\tau } \right)d\tau }}{{\sqrt{c^2 \tau ^2  - r^2}} }}} +
\]
$$
+ \frac{c}{{2\pi }}\int\limits_{\textrm{Г}_\varepsilon ^ -  \left(
x^* \right)} {H\left( {ct - r} \right)\frac{1}{r}\frac{{\partial
r}}{{\partial n\left( y \right)}}dS\left( y \right)\int\limits_{{r
\mathord{\left/ {\vphantom {r c}} \right.
\kern-\nulldelimiterspace} c}}^t {\frac{{\tau \dot u\left( {y,t -
\tau } \right)d\tau }}{\sqrt{c^2 \tau ^2  - r^2} }}}.
$$
 При $\varepsilon  \to {\rm  + 0}$ первый интеграл, в силу слабой
особенности у подынтегральной функции, стремится к интегралу по
$S$, второй - к интегралу в смысле главного значения, который
также существует, так как подынтегральная функция имеет
особенность вида $r^{ - 1} $ и содержит функцию $\frac{{\partial
r}}{{\partial n(y)}} = n_j (y)\frac{{\left( {y_j  - x_j }
\right)}}{r}$. Последняя при $y \to x$ асимптотически эквивалентна
$\frac{{\partial r}}{{\partial n(x)}} = n_j (x)\frac{{\left( {y_j
- x_j } \right)}}{r}$, которая антисимметрична относительно точки
$x$. Интегральный сомножитель  при $r \to 0$ не имеет особенности.

Рассмотрим последнее слагаемое, которое обозначим $J_\varepsilon
\left( x \right)$. На $\textrm{Г}_\varepsilon ^ -  , r =
\varepsilon , \frac{{\partial r}}{{\partial n(y)}} =  -
\frac{{\partial r}}{{\partial r}} =  - 1, dS\left( y \right) =
\varepsilon \,d\theta$ , где $\theta$ - полярный угол с вершиной в
точке $x; \theta _1 $ и $\theta _2 $- углы концов
$\textrm{Г}_\varepsilon ^ -  $, пронумерованные по порядку при
обходе контура $_\varepsilon ^{} $ в положительном направлении. С
учетом этих соотношений
$$
\begin{gathered}
J_\varepsilon  \left( x \right) = \frac{c}{{2\pi
}}\int\limits_{\theta _1 }^{\theta _2 } {\frac{\varepsilon
}{\varepsilon }d\theta } \int\limits_{{\varepsilon \mathord{\left/
{\vphantom {\varepsilon  c}} \right. \kern-\nulldelimiterspace}
c}}^t {\frac{\tau \dot u\left( {y,t - \tau } \right)d\tau
}{\sqrt{c^2 \tau ^2  - r^2}}}= \frac{{\left( {\theta _2 - \theta
_1 } \right)c}}{{2\pi }}\int\limits_{{\varepsilon \mathord{\left/
{\vphantom {\varepsilon  c}} \right. \kern-\nulldelimiterspace}
c}}^t {\frac{{\tau \dot u\left({y,t - \tau } \right)d\tau
}}{\sqrt{c^2 \tau ^2  - r^2}} }, \cr \mathop {\lim
}\limits_{\varepsilon  \to 0} \left( {\theta _2  - \theta _1 }
\right) =  - \pi ,\quad \mathop {\lim }\limits_{\varepsilon \to 0}
J_\varepsilon  \left( x \right) =  - 0,5\int\limits_0^t {\dot
u\left( {y,t - \tau } \right)d\tau  =  - 0,5u\left( {x,t}
\right)}. \end{gathered} $$
 Перенося это слагаемое в левую часть,
с учетом определения $H_D^ - (x,t) = 0,5$ для $x \in S$, получим
справедливость формулы и на границе. Поскольку слева и справа в
(\ref{(5.1)}) стоят регулярные обобщенные функции, в силу леммы
Дюбуа-Реймона ([2], с. 97) равенство, справедливое в классе
обобщенных функций, справедливо и в обычном смысле. Теорема
доказана.

Для определения неизвестной граничной функции при решении краевых
задач соотношение (\ref{(5.1)}) дает
 \textit{граничное интегральное уравнение}:
\[
\begin{gathered}
\pi c^{ - 1} u_S \left( {x,t} \right) = \\=\int\limits_{S_t \left(
x \right)} {dS\left( y \right)} \int\limits_{{r \mathord{\left/
{\vphantom {r c}} \right. \kern-\nulldelimiterspace} c}}^t
{\frac{{\partial u\left( {y,t - \tau } \right)}}{{\partial n\left(
y \right)}}\frac{{d\tau }}{\sqrt{c^2 \tau ^2  - r^2}} } +
V.P.\int\limits_{S_t \left( x \right)} {\frac{1}{r}\frac{{\partial
r}}{{\partial n\left( y \right)}}dS\left( y \right)}
\int\limits_{{r \mathord{\left/ {\vphantom {r c}} \right.
\kern-\nulldelimiterspace} c}}^t {\frac{{\tau \dot u_S \left( {y,t
- \tau } \right)d\tau }}{\sqrt{c^2 \tau ^2  - r^2}} }.
\end{gathered}
\]
В случае первой краевой задачи  левая часть уравнения и второй
интеграл справа известны, определяются граничным условием, а
первый интеграл содержит ядро со слабой особенностью, но не в
точке, как обычно для эллиптических задач, а на фронте функции
Грина. Решая его, определяем нормальную производную искомой
функции на границе, после чего формула (\ref{(5.1)}) позволяет
вычислить решение в любой точке области определения.  В случае
второй краевой задачи имеем сингулярное ГИУ для определения
неизвестных граничных искомой функции u. Решая его, определяем ее
значения на границе, после чего (5.1) определяет решение.

В случае ненулевых  начальных условий решение  задачи дает
следующая теорема.

Т е о р е м а  5.2. \textit{При $N=2$ функция решение КЗ, имеет
следующее интегральное представление}:
\[
\displaylines{ 2\pi \hat u = \int\limits_0^t {d\tau }
\int\limits_{S_\tau  \left( x \right)} {\left( {\frac{{\partial
u\left( {y,t - \tau } \right)}}{{\partial n\left( y \right)}} +
\frac{1}{r}\frac{{\partial r}}{{\partial n\left( y \right)}}\tau
\dot u\left( {y,t - \tau } \right)} \right)} \frac{{dS\left( y
\right)}}{{\sqrt {\tau ^2  - \left( {r/c} \right)^2 } }} +  \cr +
\frac{\partial }{{\partial t}}\int\limits_{S_t^ -  (x)}
{\frac{{u_0 (y)dV\left( y \right)}}{{c\sqrt {c^2 t^2  - r^2 } }}}
- ?\int\limits_{S_t^ -  (x)} {\frac{{\dot u_0 (y)dV\left( y
\right)}}{{c\sqrt {c^2 t^2  - r^2 } }}}  + c\int\limits_0^t {d\tau
} \int\limits_{S_\tau  \left( x \right)} {\frac{{G(y,t - \tau
)dV\left( y \right)}}{{\sqrt {c^2 \tau ^2  - r^2 } }}}  +  \cr -
\int\limits_{S_t (x)} {u_0 (y)\frac{{ct}}{r}} \frac{{\partial
r}}{{\partial n\left( y \right)}}\frac{{dS\left( y
\right)}}{{\sqrt {c^2 \tau ^2  - r^2 } }} }
\]

Д о к а з а т е л ь с т в о следует из теоремы 4.1. и 5.2. , если
записать свертки с начальными данными в интегральном виде. Здесь
интегралы со второго до четвертого  совпадают с формулой Пуассона
для задачи Коши. Последнее дополнительное слагаемое с начальными
данными обусловлено наличием границы области.

Используя соотношения (\ref{(5.2)}), динамический аналог формулы
Гаусса можно записать в интегральном виде.

Л е м м а  5.1. \textit{При $N=2$ динамический аналог формулы
Гаусса имеет следующий вид}:
\[
V.P.\int\limits_{S_\tau  \left( x \right)} {\frac{1}{{\sqrt {1 -
\left( {r/ct} \right)^2 } }}\frac{\partial }{{\partial n\left( y
\right)}}\left( {\ln \frac{1}{r}} \right)dS\left( y \right) +
\frac{\partial }{{c\partial t}}\int\limits_{S_t^ -  \left( x
\right)} {\frac{{dV\left( y \right)}}{{\sqrt {c^2 t^2  - r^2 } }}}
= 2\pi H_S^ - ( x )H( t ),}
\]
\textit{где интеграл в смысле главного значения  берется для
граничных точек}.

Д о к а з а т е л ь с т в о . Формулу леммы 4.1  , с учетом
(\ref{(5.2)}),(\ref{(5.3)}),  можно записать так
\begin{equation}\label{(5.4)}
- \int\limits_{S_t \left( x \right)} {\frac{{ct}}{{\sqrt {c^2 t^2
- r^2 } }}\frac{1}{r}\frac{\partial r}{\partial n\left( y
\right)}dS(y) + \frac{\partial }{{c\partial t}}\int\limits_{S_t^ -
(x)} {\frac{{\,dV(y)}}{{\sqrt {c^2 t^2 - r^2 } }}}}= 2\pi H_S^ -
(x)H(t)
\end{equation}
Покажем, что это равенство, справедливое в области регулярности ,
сохраняется и для $x \in S$, если сингулярный интеграл слева,
который содержит сильную особенность по $r$, брать в смысле
главного значения. Аналогично (\ref{(5.4)})  для области с
выколотой $\varepsilon$ -окрестностью точки $x$ получим
\[
- \int\limits_{S_\varepsilon   + \textrm{Г}_\varepsilon ^ -  }
{\frac{ct\;H\left( {ct - r} \right)}{\sqrt{c^2 \tau ^2  - r^2}}
\frac{1}{r}\frac{{\partial r}}{{\partial n\left( y
\right)}}dS\left( y \right) + \frac{\partial }{{c\partial
t}}\int\limits_{S_\varepsilon ^ - \left( x \right)}
{\frac{{H\left( {ct - r} \right)dV\left( y \right)}}{\sqrt{c^2
\tau ^2  - r^2}} }} = 0,
\]
\[
- \int\limits_{S_\varepsilon   + \textrm{Г}_\varepsilon ^ +  }
{\frac{{ct\;H\left( {ct - r} \right)}}{\sqrt{c^2 \tau ^2  - r^2}}
\frac{1}{r}\frac{{\partial r}}{{\partial n\left( y
\right)}}ds\left( y \right) + \frac{\partial }{{c\partial
t}}\int\limits_{S_\varepsilon ^ - \left( x \right) +
\textrm{Ш}_\varepsilon } {\frac{{H\left( {ct - r} \right)dv\left(
y \right)}}{\sqrt{c^2 \tau ^2  - r^2}} }}  = 2\pi.
\]
На $\textrm{Г}_\varepsilon ^ +  $  $\frac{{\partial r}}{{\partial
n(y)}} = 1$, на  $\textrm{Г}_\varepsilon ^ -  $ $\frac{{\partial
r}}{{\partial n(y)}} =  - 1$, поэтому, если оба равенства сложить
и поделить на 2, с учетом свойств симметрии подынтегральных
функций (\ref{(4.1)}), получим
\begin{equation}\label{(5.5)}
\begin{gathered}
- \int\limits_{S_\varepsilon  }\frac{ct\;H(ct - r)}{\sqrt{c^2 \tau
^2  - r^2}}\frac{1}{r}\frac{{\partial r}}{{\partial n\left( y
\right)}}dS(y) + \frac{\partial }{{c\partial t}}
{\int\limits_{S_\varepsilon ^ - ( x)} {\frac{{H( {ct -
r})dV(y)}}{\sqrt{c^2 \tau ^2  - r^2}} }}+
\int\limits_{\textrm{Ш}_\varepsilon ( x )} {\frac{H( ct - r)dV(
y)}{\sqrt{c^2 \tau ^2  - r^2}} }  = \pi
 \end{gathered}
\end{equation}
Последний интеграл вычисляется переходом к полярной системе
координат. При $\varepsilon  < ct$
\[
I_\varepsilon  \left( x,t \right) =
\int\limits_{\textrm{Ш}_\varepsilon \left( x \right)}
{\frac{{H\left( {ct - r} \right)dv\left( y \right)}}{\sqrt{c^2
\tau ^2  - r^2}} }  = \int\limits_0^{2\pi } {d\theta }
\int\limits_0^\varepsilon  {\frac{{rdr}}{\sqrt{c^2 \tau ^2  -
r^2}} }  = 2\pi (ct-\sqrt{c^2 \tau ^2  - r^2} ) \quad \Rightarrow
\]
\[
\mathop {\lim }\limits_{\varepsilon  \to 0} \frac{{\partial
I_\varepsilon  }}{{\partial t}} = 2\pi c\mathop {\lim
}\limits_{\varepsilon  \to 0} \left( {1 - \frac{ct}{\sqrt{c^2 \tau
^2  - r^2}}} \right) = 0\,\, \textrm {для} \forall t > 0.
\]
Переходя к пределу по $\varepsilon  \to 0$ в равенстве
(\ref{(5.5)}) , для $x \in S$,\, $ t>0 $ получим
\[
- V.P.\int\limits_{S_t \left( x \right)} {\frac{{ct\;}}{\sqrt{c^2
\tau ^2  - r^2}} }\frac{1}{r}\frac{{\partial r}}{{\partial n\left(
y \right)}}dS\left( y \right) + \frac{\partial }{{c\partial t}}
\int\limits_{S_t^ -  \left( x \right)} {\frac{{dV\left( y
\right)}}{\sqrt{c^2 \tau ^2  - r^2}}}  = \pi.
\]
То есть, с учетом определения $H_S^ -  (x)$, формула (\ref{(5.4)})
справедлива для любых $x$. Из этой формулы элементарными
преобразованиями получаем формулу леммы.

\textbf{6. ГИУ пространственных краевых задач $(N=3).$}

Т е о р е м а 6.1.  \textit{Для $N=3$  решение КЗ   представимо в
виде: для $x \in S^ -  $}
\[
\begin{gathered}
4\pi u\left( {x,t} \right) = \int\limits_{S_t \left( x \right)}
{\left\{ {\frac{1} {r}\frac{{\partial u\left( {y,t - {r
\mathord{\left/ {\vphantom {r c}} \right.
\kern-\nulldelimiterspace} c}} \right)}} {{\partial n\left( y
\right)}} - u\left( {y,t - {r \mathord{\left/ {\vphantom {r c}}
\right. \kern-\nulldelimiterspace} c}} \right)\frac{\partial }
{{\partial n\left( y \right)}}\frac{1} {r} + c^{ - 1} \dot u\left(
{y,t - {r \mathord{\left/ {\vphantom {r c}} \right.
\kern-\nulldelimiterspace} c}} \right)\frac{{\partial \,\ln \,r}}
{{\partial n\left( y \right)}}} \right\}} dS\left( y \right) + \\
+ c^{ - 1} \frac{\partial } {{\partial t}}\int\limits_{S_t \left(
x \right)} {u_0 \left( y \right)\frac{{\partial \,\ln \,r}}
{{\partial n\left( y \right)}}dS\left( y \right)}  + \left\{
{\int\limits_{r = ct} {\frac{{\dot u_0 \left( y \right)}} {{c^2
t}}H_S^ -  \left( y \right)dS\left( y \right)}  + \frac{\partial }
{{\partial t}}\int\limits_{r = ct} {\frac{{u_0 \left( y \right)}}
{{c^2 t}}H_S^ -  \left( y \right)dS\left( y \right)} } \right\} -  \\
- \int\limits_{S_t^ -  (x)} r^{ - 1} G(x,t - r/c)\,\textrm{для }\, x \in S; \\
\\
2\pi u\left( {x,t} \right) = \int\limits_{S_t \left( x \right)}
{\left\{ {\frac{1} {r}\frac{{\partial u\left( {y,t - {r
\mathord{\left/ {\vphantom {r c}} \right.
\kern-\nulldelimiterspace} c}} \right)}} {{\partial n\left( y
\right)}} + \frac{{\dot u\left( {y,t - {r \mathord{\left/
{\vphantom {r c}} \right. \kern-\nulldelimiterspace} c}} \right)}}
{c}\frac{{\partial \,\ln \,r}} {{\partial n\left( y \right)}}}
\right\}} dS\left( y \right)-
\\
- V.P.\int\limits_{S_t \left( x \right)} {u\left( {y,t - {r
\mathord{\left/ {\vphantom {r c}} \right.
\kern-\nulldelimiterspace} c}} \right)\frac{\partial } {{\partial
n\left( y \right)}}\frac{1} {r}} ds\left( y \right) + c^{ - 1}
\frac{\partial } {{\partial t}}\int\limits_{S_t \left( x \right)}
{u_0 \left( y \right)\frac{{\partial \,\ln \,r}} {{\partial
n\left( y \right)}}dS\left( y \right)}  + \\+ {\int\limits_{r =
ct} {\frac{{\dot u_0 \left( y \right)}} {{c^2 t}}H_S^ -  \left( y
\right)dS\left( y \right)}  +   \frac{\partial } {{\partial
t}}\int\limits_{r = ct} {\frac{{u_0 \left( y \right)}} {{c^2
t}}H_S^ -  \left( y \right)dS\left( y \right)} }  -
\int\limits_{S_t^ -  (x)} {\frac{{G(x,t - r/c)}}
{r}} dV(y) \\
\end{gathered}
\]

Д о к а з а т е л ь с т в о . При $N=3$  $\hat U(x,t)$ - простой
слой на поверхности конуса $K_t  = \left\{ {\left( {x,t}
\right):\;\left\| x \right\| = ct} \right\}$, который, в отличие
от формулы в ([4], с. 206), удобно записать в виде
\begin{equation}\label{(6.1)}
\hat U\left( {x,t} \right) =  - \frac{{\delta \left( {t - R/c}
\right)}} {{4\pi R}} , R = \sqrt {x_1^2  + x_2^2  + x_3^2 }
\end{equation}
Для любых $\varphi \left( {x,t} \right) \in D\left( {R_4 }
\right)$ (6.1) определяет линейный функционал вида:
\[
\left( {U,\varphi } \right) =  - \frac{1} {{4\pi
}}\int\limits_{R^3 } {\frac{{\varphi \left( {x,\left\| x
\right\|/c} \right)}} {{\left\| x \right\|}}dV\left( x \right)}
\]
По формулам (\ref{(3.6)}) найдем
\begin{equation}\label{(6.2)}
\hat W\left( {x,t} \right) =  - \frac{{H\left( {ct - R} \right)}}
{{4\pi R}} ,\,\, \hat H\left( {x,m,t} \right) = \frac{1} {{4\pi
}}\frac{{x_j m_j }} {{R^2 }}\left( {c^{ - 1} \delta \left( {t -
R/c} \right) + \frac{{H\left( {ct - R} \right)}} {R}} \right).
\end{equation}

Для построения интегральных представлений динамических аналогов
формул Грина и Гаусса воспользуемся следующими равенствами,
которые можно получить на основе определения операции свертки
обобщенных функций:
\[
\alpha \left( x \right)\delta \left( {t - R/c} \right) * f\left(
{x,t} \right)H\left( t \right) = H\left( t \right)\int\limits_{S_t
\left( x \right)} {\alpha \left( {x - y} \right)f\left( {y,t - {r
\mathord{\left/ {\vphantom {r c}} \right.
\kern-\nulldelimiterspace} c}} \right)dS\left( y \right)}
\]
\[
\beta \left( {x,t} \right)\delta _S \left( x \right)H\left( t
\right) * f\left( {x,t} \right)H\left( t \right) = H\left( t
\right)\int\limits_0^t {d\tau } \int\limits_{S_t (x)} {\beta
\left( {y,t - \tau } \right)f\left( {x - y,\tau } \right)dS\left(
y \right)}
\]
\[
\alpha \left( x \right)\delta \left( {t - R/c} \right) * \beta
\left( {x,t} \right)\delta _S \left( x \right)H\left( t \right) =
H\left( t \right)\int\limits_{S_t \left( x \right)} {\alpha \left(
{x - y} \right)\beta \left( {y,t - {r \mathord{\left/ {\vphantom
{r c}} \right.\,\, \kern-\nulldelimiterspace} c}}
  \right)ds\left( y \right)},
\]
\[
\alpha \left( x \right)\delta \left( {t - R/c} \right)\mathop  *
\limits_x \gamma \left( x \right)\delta _S \left( x \right) =
\frac{\partial } {{\partial t}}\int\limits_{S_t \left( x \right)}
{\alpha \left( {x - y} \right)\gamma \left( y \right)ds\left( y
\right)}
\]
\[
\alpha \left( x \right)\delta \left( {t - R/c} \right)\mathop  *
\limits_x H_S^ -  \left( x \right) = c^{ - 1} H\left( t
\right)\int\limits_{r = ct} {\alpha \left( {x - y} \right)H_S^ -
\left( y \right)ds\left( y \right)}.
\]
Вычислим свертки в (\ref{(3.7)}) с учетом этих соотношений. Первое
слагаемое равно
\[
 - \hat U\left( {x,t} \right) * \frac{{\partial u}}
{{\partial n}}\delta _S \left( x \right)H\left( t \right) =
\frac{{\delta \left( {t - R/c} \right)}} {{4\pi R}} *
\frac{{\partial u}} {{\partial n}}\delta _S \left( x
\right)H\left( t \right) = \frac{{H\left( t \right)}} {{4\pi
}}\int\limits_{S_t \left( x \right)} {\frac{1} {r}\frac{{\partial
u\left( {y,t - {r \mathord{\left/ {\vphantom {r c}} \right.
\kern-\nulldelimiterspace} c}} \right)}} {{\partial n(y)}}dS\left(
y \right)}.
\]
Второе слагаемое ($ - 4\pi \,\hat W,_j  * \dot u(x,t)n_j (x)\delta
_S \left( x \right)H\left( t \right)$) равно сумме двух сверток:
\[
\begin{gathered}
- \frac{1} {{4\pi c}}\frac{{x_j }} {{R^2 }}\delta \left( {t -
\left\| x \right\|/c} \right) * \dot u(x,t)n_j (x)\delta _S \left(
x \right)H\left( t \right) = \frac{{H(t)}}
{c}\int\limits_{S_t (x)} {r^{ - 1} \dot u(y,t - r/c)r,_j n_j (y)} dS(y) =  \\
   = c^{ - 1} H(t)\int\limits_{S_t (x)} {\dot u(y,t - r/c)} \frac{{\partial \ln r}}
{{\partial n(y)}}dS(y) - \frac{{x_j }} {{R^3 }}H\left( {ct -
\left\| x \right\|} \right)
* \dot u(x,t)n_j (x)\delta _S \left( x \right)H\left( t \right)
=\\= H(t)\int\limits_0^t {d\tau } \int\limits_{S_t (x)}
{\frac{{r,_j n_j (y)}}{{r^2 }}\dot u(y,t - \tau )} \,dS(y) =
H(t)\int\limits_{S_t (x)} {\frac{{\partial r}} {{r^2 \partial
n(y)}}dS(y)\int\limits_0^t {\dot u(y,t - \tau )d\tau } } \, =\\=
H(t)\int\limits_{S_t (x)} {u^0 (y)\frac{{\partial r^{ - 1} }}
{{\partial n(y)}}dS(y) - } H(t)\int\limits_{S_t (x)}
{u(y,t)\frac{{\partial r^{ - 1} }}
{{\partial n(y)}}dS(y)}  \\
\end{gathered}
\]
Третье слагаемое -
\[
\begin{gathered}
- 4\pi \,\hat W,_j \mathop  * \limits_x u_0 \left( x \right)n_j
\left( x \right)\delta _S \left( x \right) =  - \frac{{x_j }}
{{R^2 }}\left( {c^{ - 1} \delta \left( {t - R/c} \right) +
\frac{{H\left( {ct - R} \right)}} {R}} \right)\mathop  * \limits_x
u_0 \left( x \right)n_j \left( x \right)
\delta _S \left( x \right) =  \\
= \frac{\partial } {{c\partial t}}\int\limits_{S_t \left( x
\right)} {\frac{{y_j  - x_j }} {{r^2 }}u_0 \left( y \right)n_j
\left( y \right)dS\left( y \right)}  + \int\limits_{S_t \left( x
\right)} {\frac{{y_j  - x_j }}
{{r^3 }}u_0 \left( y \right)n_j \left( y \right)dS\left( y \right)}  =  \\
= \frac{\partial } {{c\partial t}}\int\limits_{S_t \left( x
\right)} {\frac{{\partial \ln r}} {{\partial n\left( y
\right)}}u_0 \left( y \right)dS\left( y \right)}  -
\int\limits_{S_t \left( x \right)} {\frac{{\partial r^{ - 1} }}
{{\partial n\left( y \right)}}u_0 \left( y \right)dS\left( y \right)}  \\
\end{gathered}
\]
Последние три слагаемые дают известную формулу Кирхгофа для
решения задачи Коши  [ 1,2]. Складывая эти свертки, получим
формулу теоремы.

В части, зависящей от начальных данных, в формуле теоремы
добавляется слагаемое, обусловленное наличием границы $S$. Оно
исчезает при $t > t^* (x),\, \,t^* (x) = \mathop {\max }\limits_{y
\in S} \frac{{\left\| {x - y} \right\|}}{c}$, т.к. $S_t \left( x
\right) = S$ и интеграл не зависит от $ t.  $

В формуле (\ref{(3.1)}) для $x \notin S,\, t>0$ все интегралы
существуют. Доказательство ее при $x \in S$ подобно п.5. При этом
сильную особенность имеет второе слагаемое справа. Поскольку в
этом случае на $\textrm{Г}_\varepsilon ^ -  $
\[
\mathop {\lim }\limits_{\varepsilon  \to 0}
\int\limits_{\textrm{Г}_\varepsilon ^ -  } {u\left( {y,t -
\frac{r} {c}} \right)} \frac{\partial } {{\partial n\left( y
\right)}}\frac{1} {r}dS\left( y \right) = \mathop {\lim
}\limits\textrm{Г}_{\varepsilon  \to 0} \frac{1} {{\varepsilon ^2
}}\int\limits_{_\varepsilon ^ -  } {u\left( {y,t - \frac{r} {c}}
\right)} dS\left( y \right) =  - 2\pi u\left( {x,t} \right),
\]
ясно, что для $x \in S$ (6.1) сохраняет вид, если соответствующий
сингулярный интеграл понимать в смысле главного значения и учесть
значение $H_S^ -  \left( x \right)$ на границе $S$.

Для $x \in S$ формула дает сингулярные ГИУ для решения второй
краевой задачи. Для первой краевой задачи эта же формула дает ГИУ,
но уже не сингулярное. Эти уравнения также не не являются
фредгольмовскими. Решение вопроса их разрешимости, как и в плоском
случае, предлагаю заинтересованному читателю.

Используя соотношения (3.4)-(3.6) динамический аналог формулы
Гаусса также можно записать в интегральном виде. Сформулируем
результат.

Л е м м а 6.1.  \textit{При N=3 динамический аналог формулы Гаусса
имеет следующий вид}
\[
\int\limits_{S_t \left( x \right)} {\frac{{\partial r^{ - 1} }}
{{\partial n\left( y \right)}}dS\left( y \right)}  + \frac{1}
{c}\frac{\partial } {{\partial t}}\left\{ {\int\limits_{r = ct}
{\frac{{H_S^ -  \left( y \right)}} {r}dS\left( y \right) + }
\int\limits_{S_t \left( x \right)} {\frac{{\partial \ln r}}
{{\partial n\left( y \right)}}dS\left( y \right)} } \right\} =
4\pi H_S^ -  \left( x \right)H\left( t \right)
\]
\textit{При $t > t^* \left( x \right)$ отсюда следует известная
формула Гаусса }[ 2]:
\begin{equation}\label{(6.8)}
\int\limits_S {\frac{{\partial r^{ - 1} }} {{\partial n\left( y
\right)}}} dS\left( y \right) = 4\pi H_S^ -  \left( x \right).
\end{equation}

Д о к а з а т е л ь с т в о этой формулы аналогично доказательству
леммы 5.1. Заметим, однако, что при $N=2$ из формул леммы 5.1 не
следует двухмерный аналог формулы Гаусса.

\textbf{7. Решение краевых задач для уравнения Даламбера (N=1).}

Т е о р е м а  7.1. \textit{Решение краевых задач для $N=1$ имеет
следующее интегральное представление}
\begin{equation}\label{(7.1)}
\begin{gathered}
2\hat u = c\left\{ H\left( {ct - \left| {x - a_2 } \right|}
\right)\int\limits_{{{\left| {x - a_2 } \right|}/c}}^t {u,_x
\left( {a_2 ,\tau } \right)d\tau }  - H\left( {ct - \left| {x -
a_1 } \right|} \right)\int\limits_{{{\left| {x - a_1 }
\right|}/c}}^t {u,_x \left( {a_1 ,\tau } \right)d\tau } \right\} +
\\
+ \operatorname{sgn} \left({x - a_1 } \right)H\left( {ct - \left|
{x - a_1 } \right|} \right)u\left( {a_1 ,t - \frac{{\left| {x -
a_1 } \right|}} {c}} \right)-\\ - \operatorname{sgn} \left( {x -
a_2 } \right)H\left( {ct - \left| {x - a_2 } \right|}
\right)u\left( {a_2 ,t - \frac{{\left| {x - a_2 } \right|}} {c}}
\right) +
\\
+  c^{ - 1} \int\limits_{a_1 }^{a_2 } {\dot u_0 \left( y
\right)H(ct - \left| {x - y} \right|)dy}  + u_0 \left( {x + ct}
\right)H_S^ -  \left( {x + ct} \right) + u_0 \left( {x - ct}
\right)H_S^ -  \left( {x - ct} \right).
\end{gathered}
\end{equation}

Д о к а з а т е л ь с т в о.  Обозначим $x_1  = x$.  В этом случае
([2], с.206)
\begin{equation}\label{(7.2)}
\hat U\left( {x,t} \right) =  - \frac{c} {2}H\left( {ct - \left| x
\right|} \right) ,\quad \hat U,_t  =  - \frac{c} {2}\delta \left(
{t - \left| x \right|/c} \right)
\end{equation}
\begin{equation}\label{(7.3)}
\hat W\left( {x,t} \right) =  - \frac{c} {2}H\left( {ct - \left| x
\right|} \right)\left( {ct - \left| x \right|} \right), \,\, \hat
W,_x \left( {x,t} \right) = \frac{c} {2}H\left( {ct - \left| x
\right|} \right)\operatorname{sgn} \,x ,
\end{equation}
где
\begin{equation}\label{(7.4)}
\operatorname{sgn} \,x = \left\{ \begin{gathered}
1,\quad x > 0 \hfill \\
0,\quad x = 0 \hfill \\
- 1,\;x < 0 \hfill \\
\end{gathered}  \right.
\end{equation}

В одномерном случае для построения интегрального аналога формулы
Грина непосредственно воспользоваться формулой (\ref{(3.7)})
нельзя, так как не определены некоторые из входящих в нее функций.
Можно получить аналогичную формулу, доопределяя $u$ нулем вне
заданной области и проводя дифференцирование в классе обобщенных
функций. Поступим иначе, чтобы воспользоваться формулой
(\ref{(3.5)}).

Расширим область определения $u\left( {x,t} \right)$ до полосы в
$R^2 \times R^ +  $: $\{a_1  \leqslant \;x_1 \leqslant a_2$ , $ -
\infty  < \;x_2  < \infty , t>0\}$. Тогда граница $S$ будет
состоять из двух прямых $x_1  = a_1 ,\, x_1  = a_2 $, внешние
нормали которых имеют координаты (-1,0) и (1,0) соответственно,
$\frac{{\partial u}}{{\partial n}} = n_1 \frac{{\partial
u}}{{\partial x_1 }}$ для $x \in S$], $H_S^ - \left( x \right) =
H\left( {x_1  - a_1 } \right)H\left( {a_2  - x_1 } \right)$, $n_1
\delta _S \left( x \right) = \frac{{\partial H_S^ -  }}{{\partial
x_1 }} =  - \delta \left( {x_1  - a_1 } \right) + \delta \left(
{x_1  - a_2 } \right)$.

На основе метода спуска по $x_2 $ формула (\ref{(3.7)}) (теорема
3.1) преобразуется к виду:
\[
\hat u = \hat U_2  * \frac{{\partial u}} {{\partial x}}H\left( t
\right)\left( {\delta \left( {x - a_2 } \right) - \delta \left( {x
- a_1 } \right)} \right) + \frac{{\partial \hat W_2 }} {{\partial
x}} * \dot u\left( {x,t} \right)H\left( t \right)\left( {\delta
\left( {x - a_2 } \right) - \delta \left( {x - a_1 } \right)}
\right) +
\]
\[
\hat u = \hat U_2  * \frac{{\partial u}} {{\partial x}}H\left( t
\right)\left( {\delta \left( {x - a_2 } \right) - \delta \left( {x
- a_1 } \right)} \right) + \frac{{\partial \hat W_2 }} {{\partial
x}} * \dot u\left( {x,t} \right)H\left( t \right)\left( {\delta
\left( {x - a_2 } \right) - \delta \left( {x - a_1 } \right)}
\right) +
\]
\begin{equation}\label{(7.5)}
c^{ - 2} \left( {\hat U_2 ,_t \mathop  * \limits_x u_0 \left( x
\right)H_S^ -  \left( x \right)} \right) + c^{ - 2} \hat U_2  *
\hat G
\end{equation}
Здесь уже все свертки берутся в $R^1  \times R^ +  $ либо в $R^1 $
c функцией Грина и ее первообразной при $N=2$. Свертывая по $x_2$,
получим
\[
\hat u = \hat U\left( {x - a_2 ,t} \right)\mathop  * \limits_t
\frac{{\partial u\left( {a_2 ,t} \right)}} {{\partial x}}H\left( t
\right) - \hat U\left( {x - a_1 ,t} \right)\mathop  * \limits_t
\frac{{\partial u\left( {a_1 ,t} \right)}} {{\partial x}}H\left( t
\right) +
\]
\[
\hat W,_x \left( {x - a_2 ,t} \right)\mathop  * \limits_t \dot
u\left( {a_2 ,t} \right)H\left( t \right) - \hat W,_x \left( {x -
a_1 ,t} \right)\mathop  * \limits_t \dot u\left( {a_1 ,t}
\right)H\left( t \right) + \hat W,_x \left( {x - a_2 ,t}
\right)u_0 \left( {a_2 } \right) -
\]
\begin{equation}\label{(7.6)}
\hat W,_x \left( {x - a_1 ,t} \right)u_0 \left( {a_1 } \right) -
c^{ - 2} \left( {\hat U\mathop  * \limits_x \dot u_0 \left( x
\right)H_S^ -  \left( x \right)} \right) - c^{ - 2} \hat U,_t
\mathop  * \limits_x u_0 \left( x \right)H_S^ -  \left( x \right)
+ \hat U * \hat G
\end{equation}
Подставим (\ref{(7.2)})-(\ref{(7.4)}) в (\ref{(7.6)}) и выполним
интегрирование, в результате получим (\ref{(7.1)}).

Формула выражает $u\left( {x,t} \right)$ внутри области через
начальные и граничные значения самой функции и ее первых
производных. Легко показать, что она является справедливой и для
$x = a_1 ,\,\,x = a_2 $ (с учетом определения (\ref{(7.3)})). Для
этого достаточно формулу (\ref{(7.1)}) записать для интервала $
\left( {a_1 ',a_2 } \right) = \left( {a_1  + \varepsilon ,a_2 }
\right)\, ( \left( {a_1 ,a_2  - \varepsilon } \right)$). Полагая
$x = a_1 $, делаем предельный переход по $\varepsilon  \to  + 0$:
 \[
\begin{gathered}
0 = \mathop {\lim }\limits_{\varepsilon  \to 0} c\left\{ {H\left(
{ct - d} \right)\int\limits_{\frac{{\left| {x - a_2 } \right|}}
{c}}^t {u,_x \left( {a_2 ,\tau } \right)d\tau }  - H\left( {ct -
\varepsilon } \right)\int\limits_{\frac{\varepsilon } {c}}^t {u,_x
\left( {a_1 ,\tau } \right)d\tau } } \right\} +\\+ H\left( {ct -
d} \right)u\left( {a_2 ,t - d/c} \right) + + \operatorname{sgn}
\left( { - \varepsilon } \right)H\left( {ct} \right)u\left( {a_1
,t} \right) + c^{ - 1} \int\limits_{a_1  +
 \varepsilon }^{a_2 } {\dot u_0 \left( y \right)H(ct -
 \left| {a_1  - y} \right|)dy + }  \\
+ u_0 \left( {a_1  + ct} \right)H_S^ -  \left( {a_1  + ct} \right)
+ u_0 \left( {a_1  - ct} \right)H_S^ -  \left( {a_1  - ct} \right)
=
\\
= c\left\{ {H\left( {t - d/c} \right)\int\limits_{\frac{{\left| {x
- a_2 } \right|}} {c}}^t {u,_x \left( {a_2 ,\tau } \right)d\tau }
- \int\limits_0^t {u,_x \left( {a_1 ,\tau } \right)d\tau } }
\right\} + H\left( {t - d/c} \right)u\left( {a_2 ,t - \frac{d}
{c}} \right) +  \\
- H\left( t \right)u\left( {a_1 ,t} \right) + c^{ - 1}
\int\limits_{a_1 }^{a_2 } {\dot u_0 \left( y \right)H(t - \left|
{a_1  - y} \right|/c)dy}  + u_0
\left( {a_1  + ct} \right)H_S^ -  \left( {a_1  + ct} \right) +  \\
+ u_0 \left( {a_1  - ct} \right)H_S^ -  \left( {a_1  - ct} \right) \\
\end{gathered}
\]
($d = \left| {a_1  - a_2 } \right|$). Перенося слагаемое $ -
H\left( t \right)u\left( {a_1 ,t} \right)$ в  левую часть, с
учетом значения характеристической функции на границе, получим
формулу теоремы для левой  границы. Аналогичные рассуждения
проводятся для $x = a_2 $. В результате на концах отрезка
$(a_1,a_2)$ имеем следующие соотношения для определения
неизвестных:
\[
\begin{gathered}
u(a_1 ,t) = c {H\left( {ct - d} \right)\int\limits_{d/ c}^t
{u_{,x} \left( {a_2 ,\tau } \right)d\tau }  - cH(t)\int\limits_0^t
{u_{,x} \left( {a_1 ,\tau } \right)d\tau } }  +
 H\left( {ct - d} \right)u\left( {a_2 ,t - \frac{d}
{c}} \right) +
\\
+c^{ - 1} \int\limits_{a_1 }^{a_2 } {\dot u_0 \left( y \right)H(ct
- \left| {a_1  - y} \right|)dy}  + u_0 \left( {a_1  + ct}
\right)H\left( {d - ct} \right)H\left( t \right)\quad \textrm{для}
\,\,x = a_1;
\\
u(a_2 t) = c {H\left( t \right)\int\limits_0^t {u,_x \left( {a_2
,\tau } \right)d\tau }  - cH\left( {ct - d}
\right)\int\limits_{{d/c}}^t {u,_x \left( {a_1 ,\tau }
\right)d\tau } }  + H\left( {ct - d} \right)u\left( {a_1 ,t -
\frac{d}{c}} \right) +
\\
+c^{ - 1} \int\limits_{a_1 }^{a_2 } {\dot u_0 \left( y \right)H(ct
- \left| {a_2  - y} \right|)dy}  + u_0 \left( {a_2  - ct}
\right)H\left( {d - ct} \right)H\left( {ct} \right)\quad
\textrm{для}  \,\, x = a_2.
\end{gathered}
\]

При заданных $u,_x \left( {a_k ,t} \right),\,k = 1,2$, получим два
функциональных уравнения с запаздывающим аргументом для
определения u на границе области, которые можно решать пошагово по
времени, начиная с t=0: При известных $u\left( {a_1 ,t}
\right),\,\,\,u\left( {a_2 ,t} \right)$ - это уже система
интегральных уравнений, где искомые функции   (k=1, 2) стоят под
знаком интеграла. Можно рассмотреть смешанную задачу, когда на
одном конце известна функция, а на другом ее производная. И в этом
случае эта система уравнений является разрешающей.

\textbf{Заключение.} Полученные решения КЗ ранее были опубликованы
без доказательств в [3,4]  в виде кратких сообщений (там есть
описки в знаках слагаемых). Однако автору представляется
необходимым изложение самого метода обобщенных функций  (МОФ) на
примере КЗ для волнового уравнения.  Используя его можно построить
аналогичные формулы и граничные интегральные уравнения  для
решения КЗ в пространствах большей размерности ( $N>3$).  Особенно
эффективен этот метод в  краевых задачах для систем уравнений с
частными производными с постоянными коэффициентами, когда удается
построить матрицу Грина системы. При этом несущественен тип
уравнений, он может быть и эллиптическим, как, например в [5],
параболическим или смешанного типа [6,7] .

Но особенно эффективен МОФ для решения гиперболических уравнений,
где использование аппарата классического дифференцирования весьма
затруднительно, а иногда и просто невозможно (см., например,
[8,9]). Исследование разрешимости  построенных ГИУ, которые для
задач с граничными условиями типа Неймана, являются сингулярными,
представляет самостоятельную задачу функционального анализа,
поскольку построенные уравнения не относятся к хорошо изученным
классическим. Однако заметим, что использование вычислительных
методов на основе методов граничного элемента с переходом к
дискретным аналогам ГИУ,  позволяет достаточно эффективно строить
решения подобных ГИУ [6,10]. \vspace{5mm}
\newpage
\centerline{{\textbf{Список использованной литературы}}}

1.  Петровский И.С. Лекции об уравнениях с частными производными.
М., 1961 г.

2.  Владимиров В.С. Уравнения математической физики. М., 1978, 512
с.

3. Алексеева Л.А. Гранично-интегральные уравнения начально-краевой
задачи для волнового уравнения в  . Дифференциальные уравнения,
1992 г., № 8.

4.  Алексеева Л.А. Динамические аналоги формул Грина, Гаусса для
решений волнового уравнения в
//Дифференциальные уравнения.1995, т.31, №11, cc. 1951-1953.

5. Алексеева Л.А., Саутбеков С.С.  Метод обобщенных функций при
решении стационарных краевых задач для уравнений Максвелла//
Журнал вычислительной математики и математической физики.
2000.Т.40. №4.С. 611-622.

6.  Алексеева Л.А., Купесова Б.Н. Метод обобщенных функций в
краевых задачах связанной термоэластодинамики
//  Прикладная математика и механика. 2001. T.65 . №2.С.334-345.

7.  Alekseyeva L.A.. Boundary Element Method of Boundary Value
Problems of   Elasto\-dy\-na\-mics by Stati\-ona\-ry Run\-ning
Loads
//Int. J. Engineering Analysis with Boundary Element. 1998.  N11.
P.37-44 (UK,Oxford)

8.  Алексеева Л.А. Обобщенные решения нестационарных краевых задач
для уравнений Максвелла //  Журнал вычислительной математики и
математической физики. 2002.  Т.42. №1.С.76-88.

9. Alexeyeva L.A., Zakiryanova G.K. Generalized solutions of
boundary value problems of dynamics of anisotropic  elastic
media// Journal of the Mechanical Behaviour of Materials. 2004.
№5, P.16-21.

10. Alexeyeva L.A., Dildabaev Sh.A., Zhanbyrbaev A.B., Zakiryanova
G.K. Boundary Integral Equation Method in Two and three
dimensional problems of elastodynamics// Int. J. Computational
Mechanics.1996.V.18.N2. pp. 147-157.

\end{document}